\documentstyle{amsppt}
\magnification1200
\NoBlackBoxes
\pageheight{9 true in} 
\pagewidth{6.5 true in}

\topmatter
\title 
Approximating $1$ from below using $n$ Egyptian fractions
\endtitle
\author
K. Soundararajan
\endauthor
\address
Department of Mathematics, 
University of Michigan,
Ann Arbor, MI 48109
\endaddress
\email 
ksound\@umich.edu
\endemail
\endtopmatter

\document

\noindent Define the sequence of positive integers 
$a_1=2$, $a_2=3$, $a_3=7$, $a_4=43$ and in 
general $a_{k+1}= a_1\cdots a_k +1$.  In this note we 
shall prove that if $b_1 \le \ldots \le b_k$ are natural 
numbers with $\frac{1}{b_1} +\ldots+\frac{1}{b_k} <1$ 
then 
$$
\frac{1}{b_1} +\ldots +\frac{1}{b_k} \le 
\frac{1}{a_1} +\ldots + \frac{1}{a_k},  \tag{1}
$$
and that equality in (1) holds only when $b_i=a_i$ 
for $1\le i\le k$.  Our proof is based on induction on $k$; 
the proof when $k=1$ is clear. 

Observe that $\frac{1}{a_1} +\ldots +\frac{1}{a_k} 
= 1- \frac{1}{a_1\cdots a_k}$ and that 
$\frac{1}{b_1}+\ldots +\frac{1}{b_k}$ may be expressed 
as a fraction with denominator $\le b_1 \cdots b_k$.  
Thus (1) holds strictly if $b_1\cdots b_k < a_1 \cdots a_k$ 
and we assume henceforth that $b_1\cdots b_k \ge a_1\cdots a_k$. 

Let $1\le \ell \le k$ denote the largest integer $j$ such 
that $b_j b_{j+1} \cdots b_k \ge a_{j} a_{j+1} \cdots a_{k}$.  
It follows that 
$$
b_\ell \ge a_\ell, \ \ b_{\ell} b_{\ell+1} \ge a_{\ell} a_{\ell+1}, \ \ 
b_{\ell} b_{\ell+1} b_{\ell+2} \ge a_{\ell} a_{\ell+1} a_{\ell +2}, \ \ 
\ldots . \tag{2}
$$
We shall prove that (2) implies that 
$$
\frac{1}{b_\ell} +\frac{1}{b_{\ell+1}} + \ldots +\frac{1}{b_k} 
\le \frac{1}{a_{\ell}} + \frac{1}{a_{\ell+1}} +\ldots + 
\frac{1}{a_k}, \tag{3}
$$
with strict inequality unless $b_i=a_i$ for $\ell \le i\le k$.  
But by induction hypothesis 
$$
\frac{1}{b_1} +\ldots +\frac{1}{b_{\ell-1}} \le 
\frac{1}{a_1} +\ldots + \frac{1}{a_{\ell-1}},
$$
and this is strict unless $b_i=a_i$ for $1\le i\le \ell-1$.  
Combined with (3) this proves (1).

\proclaim{Proposition}  Let $x_1\ge x_2 \ge \ldots \ge x_n >0$ 
and $y_1 \ge y_2 \ge \ldots \ge y_n >0$ be two decreasing 
sequences of $n$ positive real numbers.  Suppose that 
$y_1\cdots y_j \le x_1 \cdots x_j$ for every $1\le j\le n$.  
Then 
$$
x_1+\ldots+x_n \ge y_1+\ldots +y_n, 
$$
and the inequality is strict unless $x_i=y_i$ for all $i$.
\endproclaim 

Taking $x_1=1/a_{\ell}$, $x_2=1/a_{\ell+1}$, $\ldots$, 
$x_{k-\ell+1} =1/a_{k}$
and $y_1=1/{b_{\ell}}$, $\ldots$, $y_{k-\ell+1} =1/b_k$ 
in the Proposition we see that (2) implies (3). 

\demo{Proof of the Proposition} Set $x_{n+1} = \min(x_n, y_n) 
\frac{y_1\cdots y_n}{x_1\cdots x_n}$ and $y_{n+1} =\min (x_n,y_n)$. 
Then $x_{n+1} \le x_n$ and $y_{n+1} \le y_n$ and $x_1\ldots x_{n+1} 
=y_1 \ldots y_{n+1}$.  By scaling we may also assume that $x_{n+1} \ge 1$ 
so that all the variables are at least $1$.  We now deduce 
our Proposition from Muirhead's theorem (see Theorem 45, pages 44-48 
of [1]).  In the notation there take 
$\alpha_{i} =\log x_{i}$ and $\alpha_i^{\prime}
= \log y_i$ for $1\le i\le n+1$, and take $a_1 =e$ and 
$a_2= \ldots =a_{n+1}=1$.  The hypotheses of our Proposition 
then give the hypotheses of Muirhead's theorem ((2.18.1-3) of [1]) 
and the conclusion of Muirhead's theorem gives our desired 
inequality.

\enddemo

\enddocument